\documentclass[12pt]{amsart}%
\usepackage{amsmath,amssymb,amsthm,latexsym,graphicx}
\usepackage{amsfonts}
\usepackage{citesort}
\usepackage{amsmath}
\usepackage{amssymb}
\usepackage{graphicx}%
\providecommand{\U}[1]{\protect\rule{.1in}{.1in}}
\begin{document}
\title{Synthetic focusing in ultrasound modulated tomography}
\dedicatory{Dedicated to Jan Boman's 75th birthday}\author{Peter Kuchment}
\address{Mathematics Department\\
Texas A\&M University\\
College Station, TX 77843-3368, USA}
\email{kuchment@math.tamu.edu}
\author{Leonid Kunyansky}
\address{Mathematics Department\\
University of Arizona\\
Tucson, AZ 85721, USA}
\email{leonk@math.arizona.edu}
\maketitle

\begin{abstract}
Several hybrid tomographic methods utilizing ultrasound modulation have been
introduced lately. Success of these methods hinges on the feasibility of
focusing ultrasound waves at an arbitrary point of interest. Such a focusing,
however, is difficult to achieve in practice. We thus propose a way to avoid
the use of focused waves through the so called synthetic focusing, i.e. by the
reconstruction of the would-be response to the focused modulation from the
measurements corresponding to realistic unfocused waves. Examples of
reconstructions from simulated data are provided. This non-technical paper
describes only the general concept, while technical details will appear elsewhere.

\end{abstract}

\section*{ Introduction}

The last decade has seen the proliferation of the so called hybrid methods of
medical imaging, where different physical types of radiation are combined into
one tomographic process (e.g.,
\cite{KuKu,Wang_book,MXW_review,Tuchin,Vo-Dinh,Ammari,Scherzer}). Such a
combination allows one to alleviate deficiencies of each separate type of
waves, while combining their strengths. Thermoacoustic (and closely related
photoacoustic) tomography is an example of such a hybrid method (see
\cite{AKK,KuKu,Wang_book,MXW_review,Tuchin,Vo-Dinh,PatchSch,FR} and references therein).

In this text we discuss a very specific type of hybrid imaging, the one that
combines electrical or optical measurements with a concurrent scanning of the
object with ultrasound (US) waves. The propagating US wave slightly
modifies the local physical properties of the medium, such as its electrical
conductivity or the distribution of light scatterers. This, in turn, perturbs
the background tomographic measurements (optical or electric). These
perturbations are measured and used to reconstruct the electrical or optical
properties of the medium. Such hybrid methods hold a promise to improve the
tomographic modalities that are otherwise notorious for their instability
and/or low resolution, such as electrical impedance tomography~ \cite{Bor02}
or optical imaging~\cite{Wang_book}. For these hybrid methods we will use the
names Acousto-Electric Tomography (AET) and Ultrasound Modulated Optical
Tomography (UMOT)\footnote{Other names have also been used
\cite{Ammari,Ammari_EIT,Wang_book}.}.

It has been understood that if one could focus the ultrasound on a small spot
inside the body, knowledge of this location would have a stabilizing effect on
the reconstruction in otherwise highly unstable imaging modalities (see e.g.
\cite{AlmBK,Ammari_EIT} and Section \ref{S:EIT} below). In this text, however,
our goal is not to discuss the ultrasound modulated imaging in detail, but
rather to address the assumption of a well focused US beam. It is
known that such focusing (to the extent needed for tomography) is hard to
achieve (see, for example, the detailed discussion of this issue in
\cite{Localized}). Thus, it would be important to find a way of utilizing
instead unfocused US waves. As we will show, this can be achieved through the
synthetic focusing, i.e. by extracting the measurements corresponding to
well-focused beams, from the data obtained with unfocused waves. Under such
approach, the tomographic problem is solved in two steps (hopefully both
stable): first, synthetic focusing, and then inversion from the
\textquotedblleft focused\textquotedblright\ beam data. We will not attempt to
address the second step in this short text (see
\cite{Ammari_EIT,KuKuAET,AlmBK} where the various implementations of such a
procedure are discussed). Instead, we show that the first step (synthetic
focusing), under appropriate choice of waves, is stable and mathematically
equivalent to the reconstruction in thermoacoustic tomography (although no
thermoacoustic measurements are conducted). In some cases, this step can even
be reduced to the inversion of the Fourier transform. An example of the AET
reconstruction is also provided.

This paper is meant as a preliminary non-technical announcement, which follows
a part of the lecture given at the Jan Boman's conference ``Integral Geometry
and Tomography'' in August 2008. The technical analytic and numerical details
will be provided elsewhere.

The structure of the paper is as follows: The idea of US modulation in
tomography is summarized in Section \ref{S:modulated}. Section
\ref{S:focusing} contains the description of four different types of unfocused
US waves that one might try to use. In each case it is shown that the
synthetic focusing (within the assumed mathematical model approximation)
reduces to one of the well known analytic reconstruction procedures. Section
\ref{S:EIT} is devoted to a brief description of the AET procedure and a
numerical example. Remarks are contained in Section \ref{S:remarks}, followed
by Acknowledgments and References.

\section{Ultrasound modulation in tomography}

\label{S:modulated}

As indicated in the Introduction, the ultrasound (US) modulation in various
types of tomography is performed by sending a US wave through the object,
concurrently with the original tomographic measurement (say, optical or electric).
Varying the shape of the wave, one observes the resulting changes in the
tomographic data and tries to extract from them the quantity of interest. This
technique promises to improve significantly the stability of such modalities
as electrical impedance tomography (EIT) and optical tomography (OT).

\begin{figure}[t]
\begin{center}%
\begin{tabular}
[c]{cc}%
\includegraphics[width=1.8in,height=1.8in]{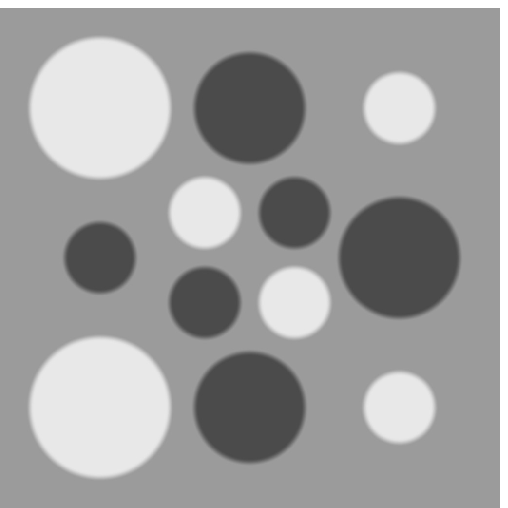} &
\includegraphics[width=1.8in,height=1.8in]{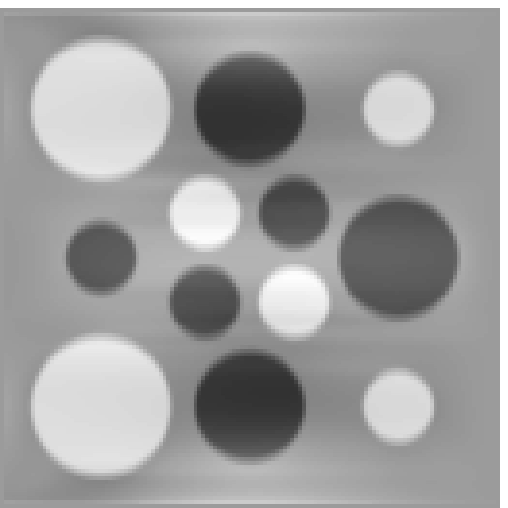}\\
&
\end{tabular}
\end{center}
\caption{Phantom (left) and its EIT reconstruction (right) using focused US
modulation.}%
\label{F:aet}%
\end{figure}

\begin{figure}[t]
\par
\begin{center}
\scalebox{0.5}{\includegraphics{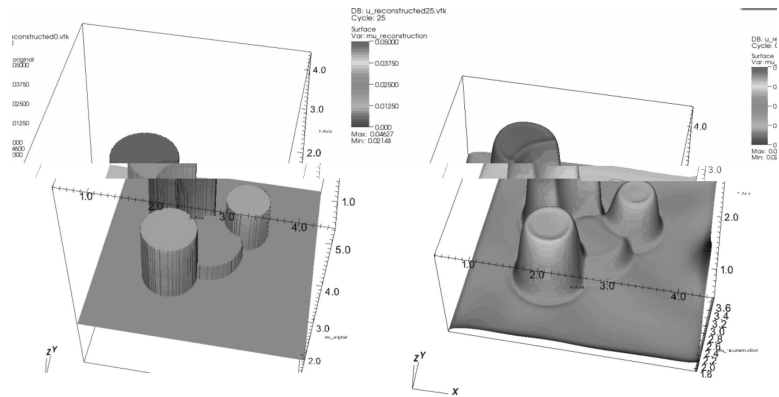}}
\end{center}
\caption{Phantom (left) and its OT reconstruction (right) using focused US
modulation.}%
\label{F:uot}%
\end{figure}

To be more specific, let us consider the EIT case \footnote{The first AET
reconstructions were probably suggested and implemented in \cite{Wang_AET}.},
where one seeks to recover the conductivity $\sigma(x)$ inside of an object
occupying the domain $\Omega$ from boundary impedance measurements. Assume for
instance that a boundary current $g(y)$ is applied and one measures the
corresponding boundary potential response $h(z)$. Here $y$ and $z$ denote
variable points on the boundary $\partial\Omega$. In EIT, the current $g$ is
varied and the potential $h$ is measured, so that the complete
Neumann-to-Dirichlet operator $\Lambda_{\sigma}$ on $\partial\Omega$ is
obtained (alternatively, by varying $h$, one recovers the Dirichlet-to-Neumann
operator). In the case of ultrasound modulation (at least in the example
discussed in Section \ref{S:EIT}), a single boundary current can be used,
rather than the whole operator $\Lambda_{\sigma}$; we will thus assume that
$g$ is fixed. Suppose now that, given any point $x\in\Omega$, we could create
a US wave that would approximate well the delta function at the location $x$.
This would create a perturbation $h_{x}(y)$ of the original boundary voltage
$h(y)$. By scanning the focusing point $x$ through the domain $\Omega$, one
obtains the set of functions $h_{x},x\in\Omega$, from which one can try to
recover the internal conductivity $\sigma$. This can actually be done, and the
reconstruction does not inherit the original high instability of EIT (see the
details in \cite{Ammari_EIT,KuKuAET}, where two different approaches are
presented). A somewhat similar situation arises in OT
\cite{Nam02,DobsonNam,AlmBK}. Stability of such reconstructions can be seen in
Figures \ref{F:aet} and \ref{F:uot} (taken from \cite{AlmBK,KuKuAET}). These
figures show mathematical phantoms (left) and their reconstructions obtained
using the simulated perfectly focused US modulations in EIT and OT
correspondingly (images on the right). These examples clearly demonstrate
stability of these hybrid modalities; the issues of stability will be
discussed in more detail, along with the reconstruction methods, in
\cite{AlmBK,KuKuAET}. 

However, as we have mentioned in the Introduction, the problem with this
technique is the difficulty of physically generating well localized US waves
(see e.g. \cite{Localized}). We thus suggest the use of a synthetic focusing
approach described in the next section.

\section{Synthetic focusing}

\label{S:focusing}

We start by mentioning that the acousto-electric and opto-acoustic effects
that are responsible for the perturbations in the electric and optical
properties of the medium due to US irradiation, are very small (e.g.,
\cite{AE2,AE3,Wang_AET,Wang_book}). Although this makes the resulting boundary
effects harder to measure, it also enables one to use the linearization of the
problem in the following sense. Let $\sigma(x)$ be the parameter to be
reconstructed, say, electrical conductivity in the EIT case. The perturbation
$\delta\sigma$ of $\sigma$ due to the applied US wave will be very small.
Hence, one can safely assume that the operator $L$ that maps the perturbation
$\delta\sigma$ into the perturbation $\delta h$ of the measured data $h$ is
linear. In fact, the operator $L$ depends on the underlying $\sigma$ and thus
will be denoted as $L_{\sigma}$.

Due to smallness of the amplitude of the applied US wave, as well as the
smallness of the acousto-electric and opto-acoustic effects, one can try also to
linearize the effect of local pressure on the perturbation $\delta\sigma$. It
turns out that such a linearization is possible in AET; it happens to be
incorrect in the optical case, where the dependence is quadratic
\cite{Wang_book}.

Let now
\begin{equation}
L_{\sigma}:\delta\sigma(x)\mapsto\delta h(y),x\in\Omega,y\in\Gamma
\subset\partial\Omega\label{E:operator}%
\end{equation}
be the linear operator that maps the distributed perturbation $\delta
\sigma(x)$ of the material parameter (say, conductivity) inside $\Omega$ to
the perturbation $\delta h(y)$ of the tomographic data measured on a part
$\Gamma$ of the boundary $\partial\Omega$ (which, in principle, could consist
of a single point, although this would not be advisable in tomography). If one
could create US waves well localized at arbitrary locations $x\in\Omega$, this
would allow one to measure the response of $L$ to delta-type inputs, and thus
would lead to a direct measurement of the kernel $l(x,y)$ of the operator
$L$:
\begin{equation}
(Lf)(y)=\int\limits_{\Omega}l(x,y)f(x)dx,x\in\Omega,y\in\Gamma.
\label{E:intoperator}%
\end{equation}
So, the use of perfectly focused US waves is equivalent to measuring the
kernel $l(x,y)$. The general idea of \textbf{synthetic focusing} is to use a
more realistic non-localized complete set $w_{\alpha}(x)$ of US waves and to
recover the kernel $l(x,y)$ from the measured responses
\[
\int\limits_{\Omega}l(x,y)w_{\alpha}(x)dx.
\]
In other words, one reconstructs numerically the would-be response of the
tomographic measurement to the perfectly focused US irradiation, which
explains our use of the name \textquotedblleft synthetic
focusing\textquotedblright.

We now address some specific versions of synthetic focusing in sub-sections
below. (The reader should be warned that we discuss here only the
reconstruction of the kernel $l(x,y)$, but not the procedure of the further
recovery of $\sigma(x)$ from the kernel.)

\subsection{Spherical waves}

\label{SS:spherical}

Suppose that one places a point US transducer at arbitrary location
$z\in\partial\Omega$ and creates a short impulse that approximates a spreading
spherical wave $p(t,z,x)=\delta(|x-z|-vt)$, where $v$ is the constant US speed
in the tissue\footnote{This certainly requires a broadband transducers, as is
the case in all other suggested synthetic focusing methods.}. Then the
response of the measurement operator $L_{\sigma}$ produces the spherical
integrals of the kernel $l(x,y)$
\begin{equation}
L_{\sigma}(p(t,z,\cdot))=\int\limits_{|x-z|=vt}l(x,y)dx,z\in\partial\Omega.
\label{E:spherical}%
\end{equation}
We thus arrive to the problem of recovering a function $l(x,y)$ of $x$ ($y$ is
now just a fixed parameter) inside domain $\Omega$ from known integrals from
$l$ over all spheres centered at the boundary $\partial\Omega$. This problem
has been studied extensively, in particular due to its relevance for the
thermoacoustic tomography (TAT); numerous results concerning the solvability,
stability, and inversion methods have been obtained (see \cite{AK,AKK,AKQ,AQ,FPR,FR,FR2,KuKu,Ku2007,PatchSch,Wang_CRC,Wang_book,MXW_review}). 
In particular,
when $\partial\Omega$ is a sphere, a variety of backprojection type formulas
have been derived. We provide below the $3D$ formula derived in \cite{MXW2}
(it was generalized to all dimensions in \cite{Ku2007}):
\begin{equation}
f(x)=\frac{1}{8\pi^{2}}\mathrm{div}\int\limits_{\partial\Omega}\mathbf{n}%
(z)\left(  \frac{1}{t}\frac{d}{dt}\frac{g(z,t)}{t}\right)  \left.
{\phantom {\rule {1pt}{8mm}}}\right\vert _{t=|y-x|}dA(z). \label{E:Kunya3D}%
\end{equation}
Here $\Omega$ is the unit sphere, $g(y,t)$ is the average of $l$ over the
sphere centered at $y\in\partial\Omega$ and of radius $t$, $\mathbf{n}(y)$ is
the outward normal vector to $\partial\Omega$ at $y$, and $dA$ is the surface
area measure.

Several different inversion formulas were derived earlier in \cite{FPR}.

In the case of non-spherical surfaces, eigenfunction expansions and time
reversal methods can be used to reconstruct $l$ (see
\cite{AQ,AK,AKK,FPR,FR,FR2} and references therein).

\subsection{Spherical monochromatic waves and plane waves in AET}

\label{SS:spherical_mono}

Since in AET the electric response of the medium is proportional to the
perturbation in the pressure, one can use for scanning not only acoustic
pulses, but also monochromatic waves. Using such waves has the advantage that
the resulting boundary measurement will oscillate with ultrasound frequency.
Hence, its time Fourier transform allows one to selectively pick up this
frequency, thus reducing the effects of the wide spectrum noise.

\subsubsection{Spherical monochromatic waves}

Suppose that the monochromatic wave $p_{\lambda}(t,z,x)$ with a time frequency
$\lambda$ is produced by a point transducer located at the point $z\in
\partial\Omega$:
\[
p_{\lambda}(t,z,x)=e^{-i\lambda t}\frac{e^{i\lambda|x-z|}}{4\pi|x-z|}.
\]
Then
\[
(Lp_{\lambda}(t,z,\cdot))(y)=\int\limits_{\Omega}l(x,y)e^{-i\lambda t}%
\frac{e^{i\lambda|x-z|}}{4\pi|x-z|}dx,x\in\Omega,y\in\Gamma,z\in\partial
\Omega,
\]
and the Fourier transform with respect to time $\widehat{Lp_{\lambda}}(z,y)$
of $Lp_{\lambda}$ will be equal to the integral of $l(x,y)$ multiplied by the
free-space Green's function $\Phi_{\lambda}(x,z)$ of the Helmholtz equation:
\begin{align*}
\widehat{Lp_{\lambda}}(z,y)  &  =\int\limits_{\Omega}l(x,y)\Phi_{\lambda
}(x,z)dx,\\
\Phi_{\lambda}(x,z)  &  =\frac{e^{i\lambda|x-z|}}{4\pi|x-z|}.
\end{align*}
As before, we assume that the point $y$ at which the electrical measurements
are conducted, is fixed. Then, if the measurements are repeated for a wide
range of frequencies $\lambda$ and for all $z\in\partial\Omega,$ and are
followed by the inverse Fourier transform (in $\lambda$) applied to
$\widehat{Lp_{\lambda}}(z,y)$, one recovers the integral of $l$ over the
spheres centered at $z$ as in (\ref{E:spherical}). Now the problem of
recovering $l(x,y)$ from values of $\widehat{Lp_{\lambda}}(z,y)$ is equivalent
to inverting its spherical mean Radon transform of $l$; \ such an inversion
was discussed in the previous section.

Alternatively, since the first step of the reconstruction formulas presented
in \cite{Ku2007} consists in computing the values of $\widehat{Lp_{\lambda}%
}(z,y)$, one can reconstruct $l(x,y)$ without explicitly reconstructing the
values of its spherical mean Radon transform first. In particular, by using the 3-D
formula from \cite{Ku2007}, one obtains the following reconstruction formula
for $l(x,y)$ from $\widehat{Lp_{\lambda}}(z,y)$:
\[
l\mathbf{(}x,y\mathbf{)=-}\frac{1}{2\pi^{2}}\mathrm{div}_{x}\int
\limits_{|z|=R}\mathbf{n(}z)h_{y}(z,|x-z|)dA(z),
\]
with%
\[
h_{y}(\mathbf{z},t)=-\frac{1}{t}\int\limits_{\mathbb{R}^{+}}\left[
\cos(\lambda t)\operatorname{Im}\left(  \widehat{Lp_{\lambda}}(z,y)\right)
-\sin(\lambda t)\operatorname{Re}\left(  \widehat{Lp_{\lambda}}(z,y)\right)
\right]  \lambda d\lambda.
\]
Here, as before, $\mathbf{n}(y)$ is the outward normal vector to
$\partial\Omega$ at $y$, and $dA$ is the surface area measure.

\subsubsection{Plane waves}

\label{SS:plane}

Spherical waves arise when the size of the transducer is small compared to the
wavelength. If the distance from the transducer to the object of interest is
much larger than the size of the object, the wave is approximately planar.
Plane waves can also be created by using a large planar transducer. In either
case the transducer should be broadband to permit generation of plane waves in
a wide range of frequencies.

Measurement done using planar acoustic waves $p_{k}(t,z,x)=e^{-i|k|t}%
\exp(ik\cdot x)$ correspond to measuring the Fourier transform of the kernel
$l(x,y)$:
\[
\widehat{Lp_{k}}(z,y)=\int\limits_{\Omega}\exp(ik\cdot x)l(x,y)dx.
\]
Synthetic focusing now reduces to the inversion of the Fourier transform.

\subsubsection{Example of AET}

\label{S:EIT}Let us illustrate by a numerical example the application of
synthetic focusing to image reconstruction in AET. We simulated numerically a
square domain, with the electrical currents equal 1 on the left and right
sides of the square and 0 on the top and the bottom. The conductivity
$\sigma(x)$ we used in our experiments was close to 1; the gray scale density
plot of the logarithm $\log\sigma(x)$ is shown in the left part of Figure
\ref{F:aet}.  In this figure the light circles correspond to the value
$\log\sigma(x)=0.05$, the dark ones represent the value of $\log
\sigma(x)=-0.05$, and the gray background depicts $\log\sigma(x)=0$. The
simulated electric potential was "measured" on the whole boundary of the
square.

\begin{figure}[t]
\begin{center}%
\begin{tabular}
[c]{cc}%
\includegraphics[width=1.8in,height=1.8in]{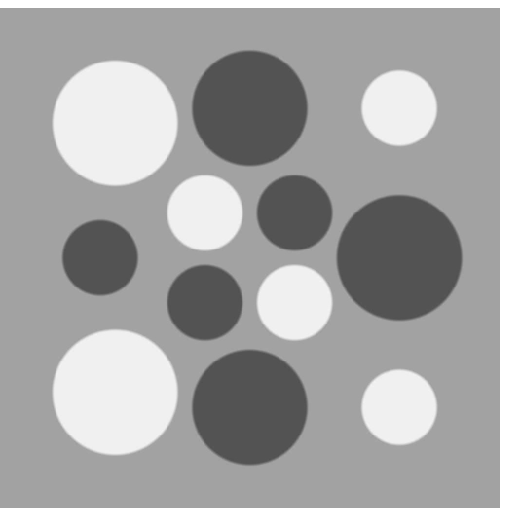} &
\includegraphics[width=1.8in,height=1.8in]{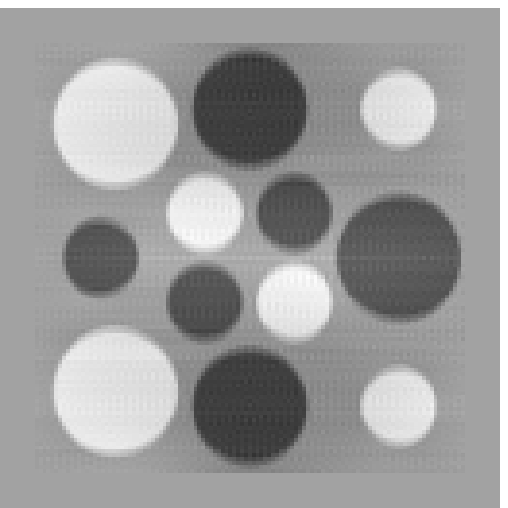}\\
&
\end{tabular}
\end{center}
\caption{AET reconstruction using spherical monochromatic acoustic waves.}%
\label{F:aetpulse}%
\end{figure}

As mentioned previously, Figure \ref{F:aet} (right) presents the result of AET
reconstruction from simulated measurements corresponding to perfectly focused
US modulation. We compare this to the results of AET reconstruction of the
same phantom, obtained by using synthetic measurements, as shown in Figure
\ref{F:aetpulse}. In this example we modeled the perturbations of electric
potential on the boundary caused by the spherical pulse waves with the centers
on a circle surrounding our square domain. The total of 300 transducers and
800 different radii of outgoing spherical pulses per transducer were
simulated. The reconstruction shown in the right half of Figure
\ref{F:aetpulse} was obtained by using synthetic focusing to obtain $l(x,y)$
and then by reconstructing conductivity $\sigma(x)$ (again, the methods for
reconstructing $\sigma$ from $l$ will be discussed elsewhere.) The
reconstructed image is clearly as good as the one obtained using the perfectly
focused US modulation (Figure \ref{F:aet}, right part).  

\subsection{Elongated focusing areas and X-ray transform}

\label{SS:elongated}

One of the UMOT reconstruction methods, suggested in \cite{LiWang_UOT2} (see
also \cite{Wang_book}), is in fact a synthetic focusing, rather than a true
reconstruction. The authors of \cite{LiWang_UOT2} employed the observation
that the focusing was non-perfect: the focal zone of the ultrasonic waves was
2 mm across and 20 mm in length, which made it a better approximation to a
segment rather than a point. One thus can conclude that by shifting and
rotating the focusing area, one can extract the X-ray transform of the kernel
$l(x,y)$ (as before, with respect to the variable $x$). Then standard X-ray
inversion formulas recover the kernel $l$.

\section{Remarks}

\label{S:remarks}

\begin{enumerate}
\item As it has already been mentioned, exact analytic and numerical details
of the reconstruction procedures will be provided elsewhere.

\item Practical implementation of the suggested synthetic focusing procedures
will have to include more precise information about possibilities and
imperfections of the currently available transducers.

\item The UMOT situation differs from AET and will probably not allow usage of
monochromatic US modulation. However, the modulation by spherical pulse waves
should still be possible. These issues will be addressed elsewhere. (As it has
been already mentioned, the synthetic focusing using X-ray transforms was
implemented for UMOT in \cite{LiWang_UOT2,Wang_book}.)
\end{enumerate}

\section*{Acknowledgments}

The work of the first author was partially supported by the NSF DMS grant
0604778 and by the KAUST grant KUS-CI-016-04. The work of the second author
was partially supported by the DOE grant DE-FG02-03ER25577. The authors
express their gratitude to NSF, DOE, and KAUST for the support. The authors
also thank M.~Allmaras, G.~Bal, W.~Bangerth, J.~Schotland, L.-H.~Wang, and
Y.~Xu for useful information and discussions. We are grateful to M.~Allmaras
and W.~Bangerth for allowing us to use Fig. \ref{F:uot}.


\begin{thebibliography}{99}
\bibitem {FR2}D.~Finch and Rakesh, Recovering a function from its spherical
mean values in two and three dimensions. To appear in \cite{Wang_CRC}.
\end{thebibliography}
\end{document}